\input amstex
\documentstyle{amsppt}
\input epsf
\input psfrag

\define\RR{\Bbb R}
\define\CC{\Bbb C}
\define\SS{\Bbb S}

\topmatter
\title Polynomial Roots and Open Mappings
\endtitle
\rightheadtext {}
\author Jon A. Sjogren
\endauthor
\affil Faculty of Graduate Studies, Towson University
\endaffil
\address Towson, Maryland
\endaddress
\date 1 June 2015
\enddate
\endtopmatter
\document

\abovedisplayskip=13pt
\belowdisplayskip=13pt

\head Plan of Argument \endhead
We examine a class of proofs to the Fundamental Theorem of Algebra that relate to partial open mappings of the complex plane.
These proofs use the ``open'' property  of a complex polynomial, at points of its domain. Based on proofs by F.S. Cater and D. Reem, we know using only elementary analysis that a non-constant polynomial is in fact open.
This fact, which is easily derived using the non-elementary tools of complex analysis, may be more technical than necessary for proving FTA. In any case,  combined with the Principle of S. Reich that for a polynomial $P(z)$, the image set $P(\CC)$ is also {\it closed} (the mapping $P(z)$ is proper), the proof is already finished!

Another approach is to exploit the openness that holds away from critical  points and critical values. What is needed follows from  the Inverse and Implicit Function Theorems. One of the simplest and best of all proofs of FTA is due to \cite{Wolfenstein}, who shows that critical points and values are easily dealt with. A related proof given below is reminiscent of the Argand-Cauchy-Littlewood method, where one shows that if $P(z_0) \ne 0$, some $z_1 \in \CC$ can always be found so that $|P(z_1)| < |P(z_0)|$. It turns out that this proof was already sketched out in Smale's survey article \cite{Smale}, where the author tacitly makes use of Reich's Principle, so it could be called the Reich-Smale Proof.

The proof of FTA from J. Milnor's published notes ``Topology
from the Differentiable Viewpoint'', \cite{TFDV}, is examined next as it also deals with topological properties of the given polynomial mapping $P(z)$. The fact that $P(z)$ is proper (has the ``propriety'' property) is used to compactify the mapping, making it possible to use a ``Pre-Image'' result from differential topology to establish that the mapping, in the non-constant case, is surjective. This argument, which uses the ``locally constant'' nature of a counting function for the pre-image points of the compactified polynomial, goes forward because of propriety. In fact it is Reich's Principle that allows  $P$ to be extended to a compact manifold (the two-sphere $\SS^2$).

The final word on giving a modern cast to Gauss's first proof (Thesis 1799) appears in \cite{Gersten-Stallings}. The blend of differential topology with the geometric theory of {\it free groups} may have resulted in the paper not being much quoted, except by Martin, Savitt and  Singer, referred to as \cite{MSS}, who study the combinatorics of
harmonic functions and their graphs. These authors come up with a new proof of FTA along the lines of Gauss I. This proof has elegant features but does not dwell on the rigor of the topological argument. It takes as second nature several observations on plane curves from the algebraic side that were overlooked \linebreak in \cite{Gersten-Stallings}.

It should be enough for grounding in the Gauss I proof, to read \cite{Ostrowski}, \cite{Gersten-Stallings}, \cite{MSS} and the present report. From the Ostrowski paper, only the first sections on ``the exterior of the large circle'' are used here, or much quoted by \cite{Smale} or  by others.  A scholar with fortitude can attempt the Master's original \cite{Gauss}. In addition, references such as the text \cite{Guillemin \& Pollack}, a source of multivariate methods such as  \cite{C.H. Edwards}, and a primer on plane curves \cite{G. Fischer} will be helpful.

As pointed out by S. Smale in his survey, the assumptions without proof made by Gauss about algebraic curves were not dealt with until the 1920's, long after the development of much of ``higher'' algebraic geometry. Our current understanding of point sets, contraction mappings etc., no doubt makes the task so brilliantly attacked by Gauss more tractable. We hope that further improvement along indicated lines will be forthcoming.  The author thanks Prof. Chao Lu and colleagues at Towson University for sharing related algorithmic work.  The author also is grateful to Dr. Daniel Reem of IMPA for his keen interest in the topic, and for suggesting a number of specific improvements to the presentation.

\head The Principle of S. Reich  \endhead

A straightforward way to prove the Fundamental Theorem of Algebra is to observe that if $P(z)$ is a polynomial of degree $\geq 1$, the image set $P(\CC)$ is both open and closed in $\CC_w$ (the ``target plane''). Since the latter ``space'' is also connected, we have $P(\CC_z) = \CC_w$, where $\CC_z$ is the ``source'' complex plane, and certainly $O_w = 0 + i 0 \in \CC_w$ belongs to the image.

In fact {\it any} complex polynomial gives a {\it closed} mapping \cite{S. Reich}. In the case of a constant $P$ (deg $P = 0$), this result is clear. For  deg $P(z) = n \geq 1$ we know that \cite{Hille} for some $R > 0$, any $z$ with $|z| \geq R$ yields $|P(z)| \geq \frac{1}{2} |z|^n$, hence $\{|z_i|\} \to \infty$ implies $\{|P(z_j)|\} \to \infty$. This means by definition that $P : \CC_z \to \CC_w$ is {\it proper} (the inverse image of a compact set $K \subset \CC_w$ is always compact).

We now state a usable form of the Inverse Function Theorem. An ``analytic'' function defined on an open set $U$ is one that has a convergent power series on $U$.

\bigskip
{\bf Proposition.} {\it Suppose $f: U \subset \CC_z \to \CC_w$ is analytic on $U$ with $f'(p) \ne 0$, $p \in U$. Then there is an open neighborhood $V$ of $v = f(p)$ and an analytic function $g = V \to U$ such that $z \in g(V)$ implies that $z = g \circ f(z)$ and $w \in V$ implies $w = f \circ g (w)$.} \qed

\bigskip
Thus at a {\it regular point} of $f$, where the derivative does not vanish, an  analytic inverse can be found on some neighborhood of the image value. The inverse function $g$ maps $V$ injectively onto an open set of $\CC_z$. We will subsequently derive this Proposition from the  {\it Implicit} Function Theorem.

One of the differences between $\RR$ and $\CC$ is that while polynomials defined on either field are continuous, proper, and closed, a polynomial with an extremum at $x_0 \in \RR$ is {\it  not open} on a small neighborhood of $x_0$.

\newpage
The use of contour integrals allows for a beautiful explicit formula for the inverse such as
$$g(v) = \frac{1}{2 \pi i} \oint_{\gamma} \frac{t f' (t)}{f(t)-v} dt$$
for a simple contour $\gamma$ contained in $U$ that holds $v$ in its interior.
In any case we have $f(p)$ contained in a $\CC_w$-neighborhood $V$ of {\it image values} meaning that $f$ is open away from the set $B \subset \CC_z$ of {\it singular points}, where the derivative {\it vanishes}. We are ready for

\bigskip {\bf Fundamental Theorem of Algebra.} (\cite{Wolfenstein, 1967}).

{\it Proof.} Let $P(z)$ have degree $n \ge 0$. Let $S = P(\CC_z)$ which is closed in $\CC_w$ by Reich's Principle so $\CC_w -S$ is {\it open}. As long as $P$ is not constant, $B := \{z: P'(z) = 0\}$ is {\it finite} in $\CC_z$ since $P'(z)$ has degree $n-1$. Hence  $T = P(B)$ is finite and $\CC_w -T$ is connected. (There is a polygonal arc connecting $w_1, w_2 \in \CC_w -T$ within that space, see \cite{Dugundji, V.2.2})

Every $w_0 \in S-T$ satisfies $w_0 =P(z_0)$ for some $z_0 \in \CC_z -B$, hence is a ``regular value'' for $P$. The Inverse Function Theorem now asserts that some neighborhood of $w_0$ maps analytically by a function locally the inverse of $P(z)$ onto a neighborhood of $z_0$. Therefore, $S-T$ is open in $\CC_w$. Writing
$$\CC_w -T = (\CC_w-S) \cup (S-T),$$
we are faced with a disjoint union of open sets. The left-hand side is connected, so if the image $S$ does not fill up $\CC_w$, we must have $S-T = \emptyset$.
But $S$ is the continuous image of  connected $\CC_z$, so is connected, and $T$ is discrete, so $S$ itself must be a one-point set $\{w_0\}$. Hence $P(z)$ must be a constant function (of degree $0$), otherwise $P(z)$ is surjective and certainly has a root. \qed

Before we move on to methods needing deeper concepts from topology and polynomial algebra, we consider a new proof combining elements of the Wolfenstein proof with another one found in \cite{Thompson}. The visualization of this new proof, which could be called the Reich-Smale approach, may appeal to some researchers.

{\it Proof of} {\bf Reich-Smale FTA.} Consider a bounded neighborhood $E$ in a $45^{\circ}$ sector  of $\CC_z$. As $z$ ranges over $E$, the image values $\{P(z)\}$ can be made to range within a $45^{\circ}$ sector of $\CC_w$ by shrinking $E$. We assume that $0 \in \CC_w$ is {\it not} in the image. Also we are assured that not all of $E$ maps to a particular radial line (ray to the origin). For one may pick a point $z_0 \in E$ where $P'(z_0) \ne 0$, else the derivative would vanish on an open set and the polynomial would be degenerate (constant). Furthermore let                    $E_0 \subset E$ be open and contain $z_0$. Then the Inverse Function Theorem implies that some open set of $\CC_w$ containing $w_0 = P(z_0)$ is the analytic image under $P$ of a subset of $E_0$. In particular, $P$ is continuous, injective and surjective from the subset of $E_0$ to its image. Thus a line segment such as part of a radial ray cannot contain the image of $E_0$.  Note how it is important to keep the image from submerging onto a radial, whereas in the Cater proof to follow, about the openness of a polynomial function, one places image points {\it onto} a radial ray.

Now we take $z_1, z_2$ close enough in $E$, and connect them with a line segment $l \subset \CC_z$, with $P(l)$ not on any radial, hence $w_1 = P(z_1)$, $w_2 = P(z_2)$ have diverse arguments (complex phases differ). By compactness, the set $P(l)$ realizes its  maximum modulus at a value $w_{*} = P(z_*)$. Now by the propriety of $P$, we can choose $R$ large enough so that $|z| > R$ implies that $|P(z)| > |w_*|$, and let $A_R := \{z : |z| \leq R\}$.
Then we know that between $\theta_1 = \text{arg}\,w_1$ and $\theta_2 = \text{arg}\,w_2$, every $\theta$ with $\theta_1 \leq \theta \leq \theta_2$ leads to  $\rho e^{i \theta} = P(z)$ for some $z \in A_R$. We make take without loss of generality the angles $0 \leq \theta_1 < \theta_2 < 2 \pi$ as lying within some $45^{\circ}$ circular sector. A real quantity depending on $\theta$ is
$$\rho_{\theta} = \inf \left\{\rho: w = \rho e^{i \theta} = P(z)\quad \text{for some}\; z \in A_R\right\}.$$
The modulus $\rho_{\theta}$ is attained by $P(z)$ on $A_R$ so we can define
$$w_{\theta} = \rho_{\theta} e^{i \theta}.$$
\indent For the continuum of allowed  values of the parameter $\theta$, not all $w_{\theta}$ can be critical values of $P$! Thus we pick $\hat{\theta} \in [\theta_1, \theta_2]$ where $w_{\hat{\theta}}$ is a regular value.  See [Figure A].

\bigskip
\centerline{\psfrag{A}{$A_R$}
\psfrag{S}{Sector}
\psfrag{E}{$E$}
\psfrag{R}{$R$}
\psfrag{l}{$l$}
\psfrag{1}{$z_1$}
\psfrag{2}{$z_2$}
\psfrag{O}{$O$}
\epsfbox{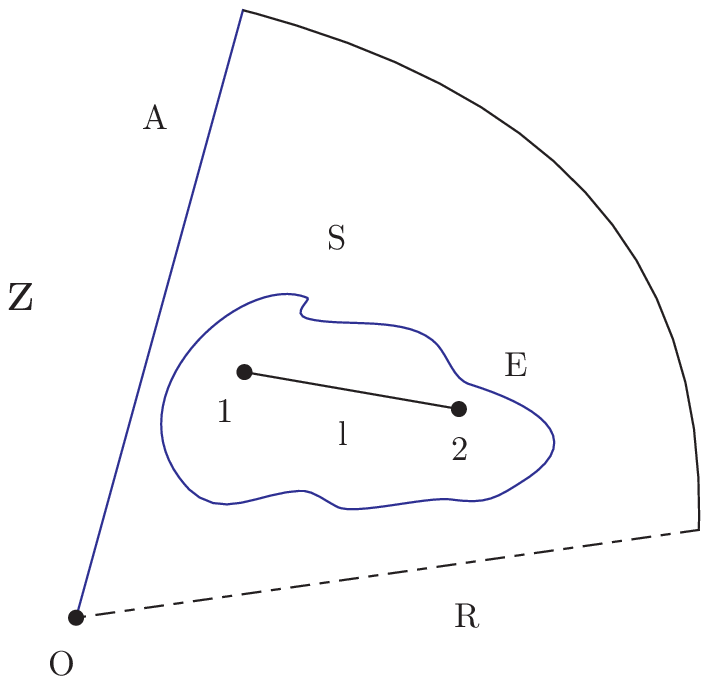}}

\vfill

\centerline{\psfrag{O}{$O$}
\psfrag{1}{$\theta$}
\psfrag{2}{$\theta'$}
\psfrag{3}{$w_{\theta'}$}
\psfrag{4}{$w_{\theta}$}
\psfrag{5}{$w_2$}
\psfrag{6}{$\hat{\theta}$}
\psfrag{7}{$w_*$}
\psfrag{8}{$P(l)$}
\psfrag{9}{$w_1$}
\epsfbox{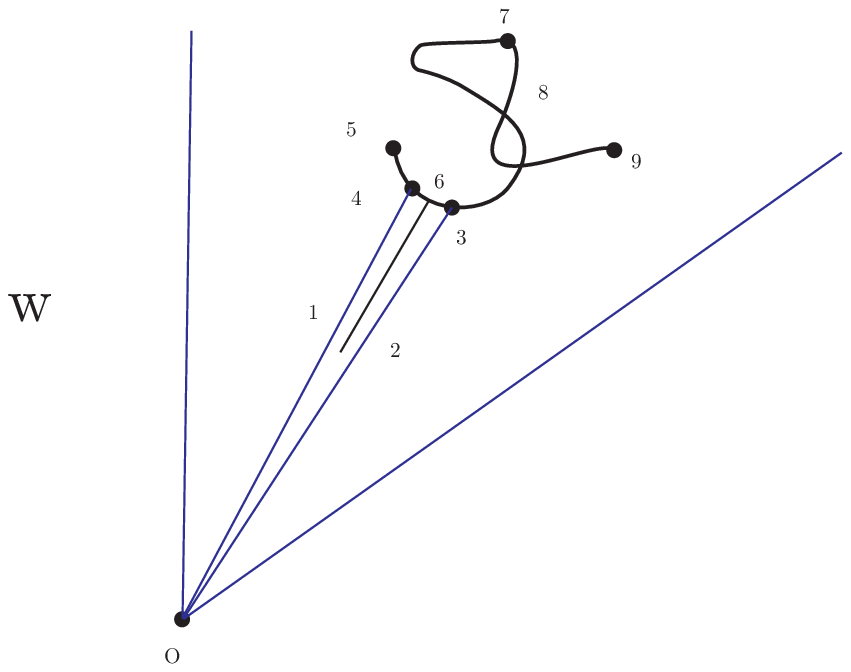}}

\centerline{Figure A}

\newpage

On the one hand, by minimality, there is no $0 < \rho <  \rho_{\hat{\theta}}$ such that there exists $z' \in A_R$ with $P(z') = \rho e^{i \hat{\theta}}$. On the other hand, $w_{\hat{\theta}}$ as a regular image point under $P$, is interior to an open set of other regular values, each having a pre-image. Thus there is a value  $P(z'')$, $z'' \in A_R$ having the argument $\hat{\theta}$, which is closer to the origin $O_w$ than is $P(\hat{z})$.

The point $z''$ must lie in $A_R$, since the  complement of $A_R$ in $\CC_z$, call it $B_R$, maps entirely to values $w=P(z)$ satisfying
$$|w| > |w_*| \geq |P(\hat{z})| > |P(z'')|\,\,.$$
The construction of $z''$ gives a contradiction which shows that $P: \CC_z \to \CC_w$ must have a root, and in fact is surjective.\qed

The present ``original'' proof of FTA just given can be seen as simplifying the \cite{Thompson} and the Milnor \cite{TFDV} proofs, and is somewhat less abstract than the Wolfenstein proof. It retains the essence of the classic Argand-Cauchy-Littlewood method, (see \cite{Littlewood}),  but is is really the same as a proof sketched out  in the ``computational'' survey of \cite{Smale}. The  author Prof. Smale uses ``propriety'' or the ``closed property'' implicitly, and did not banner the result as a theorem, so we propose to call what we just detailed, the ``Reich-Smale proof'' of FTA.

Again, the main point is the predominance of regular points and regular values. Although an arbitrary complex polynomial turns out to be an open mapping $\CC_z \to \CC_w$, this is a more obvious (local) fact when the complex  derivative is non-zero. In that case the  Inverse Function Theorem shows that a ``compact'' image set avoiding the origin must actually be ``open'' and hence include it. This contradiction shows that $A_R$ must contain a root of $P$, or else be constant.

\head Remarks on Milnor's Proof \endhead

In his ``iconic'' set of lecture notes of 1965, ``Topology from the Differentiable Viewpoint'' \cite{TFDV}, J. Milnor gave a new proof of FTA that received wide attention. The Fundamental Theorem of Algebra was displayed as an exercise in the category of smooth compact manifolds and their mappings.
Milnor resorts to the artifice of compactifying complex planes $\CC_z$ and $\CC_w$, which necessitates conjugation by stereographic projections etc. But working with compact spaces is consistent with the theme of the book \cite{TFDV}.

Consider the smooth mapping $f: \SS^2 \to \SS^2$, derived from the original polynomial $P: \CC_z \to \CC_w$. One observes that the sets of critical points $\{\kappa\} \subset \CC_z$ and critical values  $\{\tau\}\subset \CC_w$ are discrete and finite, so their complements are connected, similar to as comes up in Wolfenstein's proof. An interesting aspect of the  proof is the author's verification that $f$ is smooth at the North Pole of $\CC_z$, where in fact $f(\text{North}_z) = \text{North}_w$. This amounts to nothing less than the Reich Principle ($f$ is a closed mapping), so essential to an FTA proof of the type we are considering.

 We quote freely from \cite{TFDV}. Given $f: M \to N$ smooth,

``for ... a regular value $y$, we define $\#f^{-1}(y)$ to be the number of points in $f^{-1}(y)$.''

``The first observation to be made about $\#f^{-1}(y)$ is that it is locally constant as a function of $y$ (running through regular values!). I.e., there is a neighborhood $V \subset N$ such that $y' \in V$ implies $\# f^{-1}(y') = \# f^{-1}(y)$.''

A brief demonstration of this last quotation starts with ``let $x_1, \dotsc, x_k$ be the points of $f^{-1}(y)$, and  choose pairwise disjoint neighborhoods $U_1, \dotsc, U_k$ of these...'' From this results the invariance of the integer $k$ as $y$ varies over an open set. In the case $f^{-1}(y) = \emptyset$, the omitted proof would be that since $M$ was assumed compact, $f$ must be proper, so its image is closed in $N$.  Hence the set of ``non-image'' regular values $\{ y \}$ is also open.

A function such as $F(x, y) = \left(e^{x}\cos y, e^x \sin y\right)$ cannot be extended continuously to the Riemann sphere $\SS^2$. Indeed $(0,0) \in \RR^2_w$ is an isolated regular value with no pre-image, but every other value in $\CC_w$ {\it does} have some pre-image. The derivative matrix always has full rank, since the derivative $e^z$ never vanishes. The fact that non-image values form an open set is critical to Milnor's proof of FTA. A function such as $e^z$ is not a proper mapping from $\CC_z \to \CC_w$. On the other hand, the set $\{w: \forall z, P(z) \ne w\}$ is the complement of the image set $S = P(\CC_z)$ hence is open, since $P$ is continuous and proper, hence  closed. Thus for a polynomial $P(z)$, and $f$ derived from it, $\#P^{-1}(w)$ and $\#f^{-1}(w)$ {\it can} be defined and is locally constant.

Although there are by now a number of proofs of FTA referring to open mappings and critical points, Milnor's book shows the relevance of modern differential topology. The Pre-Image Theorem is central to the treatment of the Gauss Proof (Thesis, Univ. Helmstedt 1799) as in \cite{Gersten-Stallings}; instead we emphasize the related but more basic Implicit Function Theorem.

\head The Complex Polynomial as an Open Mapping \endhead

 We review the argument that  a non-constant complex polynomial maps every {\it open} planar set onto another open set. The simplifications from modern proofs that cover any {\it analytic} function are not substantial ---one can apply them to the case of a ``finite power series'' or polynomial. The desired result, which immediately yields that $P(\CC_z)$ is open, leads to the Fundamental Theorem of Algebra that $P(\CC_z) = \CC_w$ for deg $P(z) \geq 1$, since this space is also {\it closed} by Reich's Principle, and $\CC_w$ is connected.
\bigskip
{\bf Theorem (Complex Polynomials)}.  {\it If $f(z) = z^n + a_{n-1} z^{n-1} + \cdots + a_j z^j$, $n\geq 1$, $j \geq 0$, $a_i \in\CC$ for $i = j$, $j + 1, \dotsc, k-1$, and $a_j \ne 0$.

Then there exists $\delta > 0$ real such that $|w| < \delta$ implies that $w \in f(\CC_z)$; hence  $f$ is open at $O_z$. Thus our situation is that the polynomial $f$ is of  degree $n$ (we may write $a_n = 1$ if called for), is not constant on any neighborhood, and satisfies $f(0) = 0$. By translating $f(z)$ we see at once that $f(\CC_z) \subset \CC_w$ is also an open set.}

\bigskip
{\it Proof of Theorem.} Take $r > 0$ so small that
$$
r^n + \sum_{i = j+1}^{n-1} |a_i| r^i < |a_j| r^j \quad (a_j \ne 0). \tag{1}
$$
The inequality (1) will also hold for smaller $r'$, $0 < r'< r$. Next, consider the closed $r$-disc $D \subset \CC_z$, $D = \{z: |z| \leq r\}$ with $B = \partial D$, the complex numbers of modulus $r$. Since $f$ is continuous, $f(B)$ is compact and $\inf |f(b)|$ is realized at some $b_0 \in B$ and we let $d = |f(b_0)|$ with  $\delta = d/2$.

Given $t \in \CC_w$ satisfying $|t| < \delta = d/2$, it is sufficient for the conclusion of the Theorem to show that $t \in f(D)$. Let us assume otherwise, by compactness we can realize $\inf \{|t -f(z)|: z \in D\}$ by the choice of  some $v \in D$, letting $q = t-f(z)$, noting that $|q| > 0$.
But $v$ is {\it not} in $B$, since the value $t$  is closer to $f(O_z) = O_w$ than to $f(B)$, given $\inf \{|t -f(b)|: b \in B\} \geq d/2$, since $t$ was chosen closer in to $O_w$ than half the radius $d = \inf \{|w|: w \in f(B)\}$.  See [Figure B].

\bigskip
\centerline{
\psfrag{t}{$t$}
\psfrag{1}{$\delta$}
\psfrag{2}{$O_w$}
\psfrag{3}{$f(b_0)$}
\psfrag{4}{$f(B)$}
\psfrag{5}{$t-f(v)$}
\psfrag{6}{$b_sp^s$}
\psfrag{7}{$b_s$}
\psfrag{8}{$O_w$}
\psfrag{9}{$p = \beta e^{i \psi}$}
\epsfbox{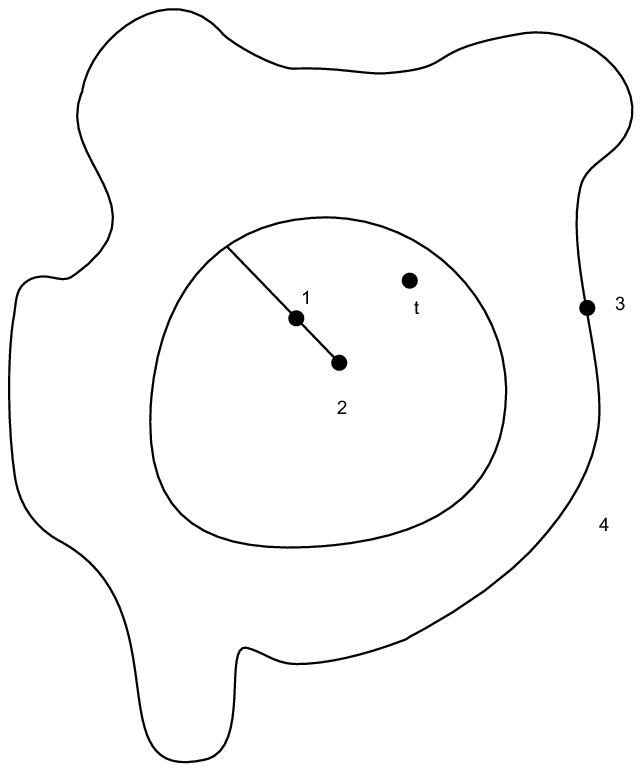}}
$$\text{Figure B}$$

\bigskip
In any case $f(b)$ never takes the value $O_w$ when $r$ is chosen small as above. For that would entail, for $|z| = r$, that
$$\left|a_j z^j\right| = |a_j| r^j = \left|z^n \sum_{i = j+1}^{n-1} a_i z^i\right| \leq r^n + \sum _{i = j+1}^{n-1} |a_i| r^i$$
which violates our postulated inequality (1).

We recapitulate the situation regarding  points and values. The value $t \in \CC_w$ was chosen closer to $f(O_z) = O_w$ than to any $f(b),\, b \in B$. We find $v \in D$ that makes $|t -f(v)|$ minimal. Hence $v$ cannot belong to $B = \partial D$ since $O_z \in D$ and $|t - f(O_z)|$ would be smaller than this minimum.
In short, $v \in D \backslash B$ so $|v| < r$.

We are now in a position to recalibrate  the function $f$ with $v$ as the new base point of a Taylor series. In other words we obtain $f(v+h) = h^n + b_{n-1} h^{n-1} + \cdots + b_s h^s$ where $s \geq 0$, $b_n \geq 1$ and $b_s \ne 0$. We note that the degree of $f$ in the new variable did not change.

\newpage
Since $|v| < r$, we may choose $\beta \in \RR$ so that all the following hold:
$$\align
\text{\it i)} & \;0 < \beta < r -|v| \\
\text{\it ii)} & \;\sum_{k = s+1}^n |b_n| \beta^{k-s} < |b_s|, \;\text{and since}\; q = t - f(z) \ne 0.\\\vspace{5pt}
\text{\it iii)} & \;|b_s| \beta^s < |q|
\endalign
$$
\indent
Now express in polar form:
$$\align
b_s &= |b_s| e^{i \theta} \quad 0 \leq \theta < 2 \pi \\
q &= |q| e^{i \varphi} \quad  0 \leq \varphi < 2 \pi \ ,
\endalign$$
and choose $0 \leq \varphi < 2\pi $ such that
$$s \cdot \psi + \theta = \varphi \mod (2\pi) \, .$$

\centerline{
\psfrag{t}{$t$}
\psfrag{1}{$\delta$}
\psfrag{2}{$O_w$}
\psfrag{3}{$f(b_0)$}
\psfrag{4}{$f(B)$}
\psfrag{5}{$t-f(v)$}
\psfrag{6}{$b_sp^s$}
\psfrag{7}{$b_s$}
\psfrag{8}{$O_w$}
\psfrag{9}{$p = \beta e^{i \psi}$}
\epsfbox{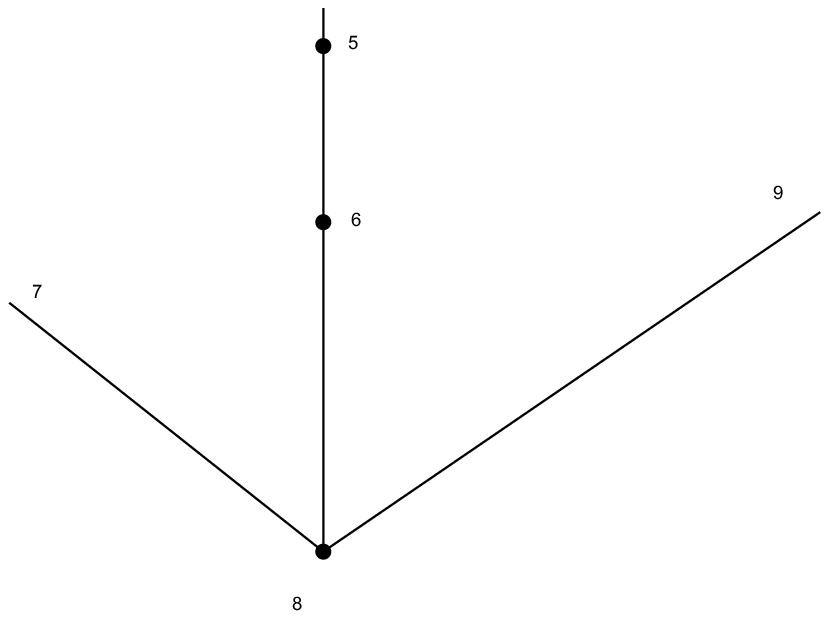}}
$$\text{Figure C}$$

Take $p = \beta e^{i \psi}$, so that $b_s p^s$ and $q = t - f(v)$ have the same argument $\mod (2\pi)$. See [Figure C]. Using {\it iii)} above, $|b_s p^s| < |q|$. Therefore since these values lie on a common radial,
$$|q| -|b_s p^s| = |t -f(v) -b_sp^s|$$
and from {\it ii)},
$$\left|\sum_{k = s+1}^n b_k p^k\right| < |b_sp^s|\, ,$$
so we finally obtain
$$\align
|t - f(v+p)| &= \left|t - f(v) -b_s p^s - \sum_{k = s+1}^n b_k p^k\right| \leq
 |t -f(v) - b_sp^s| + \left|\sum_{k = s+1}^n b_k p^k\right| = \\\vspace{2.5pt}
 & = |t-f(v)| -|b_sp^s|  + \left|\sum_{k = s+1}^n b_k p^k\right| < |t - f(v)| \ .
\endalign
$$

We must additionally understand why $v+p \in D$ holds true. But
$$|v+p| \leq |v| + |p| = |v| + \beta < r$$
by $i)$. Thus $v+p \in D$, but $|t-f(v+p)|$ being strictly smaller than $|t-f(v)|$, defined as the infimum of all  $|t-f(z)|$, $z \in D$, gives a contradiction. Hence $t \in f(D)$ after all and $O_w = f(O_z)$ is contained in a neighborhood $U(O_w)$ entirely in the image $f(D) \subset f(\CC_z)$. Translating  the polynomial in the image plane as necessary we recover

\medskip
{\bf Theorem.} {\it When $f: \CC_z \to \CC_w$ is a polynomial function of degree $n \geq 1$, if $U$ is open in $\CC_z$ then $f(U)$ is open in $\CC_w$. In particular $f(\CC_z)$ is always an open set in the complex metric topology.}\qed

\head Intersection Geometry of Plane Curves \endhead

For a real (plane) algebraic ``curve'' defined by $F(x, y)= 0$ where $F$ is a polynomial in two variables, it is often desired to exhibit some part of the locus as a smooth curve in the sense of analytic geometry. In particular, Gauss's First Proof \cite{Gauss, 1799} exhibits solutions to $f(z) = 0$ as the common points of two smooth curves in a planar domain $D \subset \RR^2$ (each with several components).
Gauss's speculation on the topological nature of these real curves was evidently premature, the real numbers not yet having  been precisely defined. Criticism has continued  until the present day with the ``completion'' work of Ostrowski clarifying some but not all of the obscurities in Gauss's arguments.
In fact Uspensky, who probably knew the work of Ostrowski as well as anyone, gives only a summary of it in his book. Smale in his survey lauds the Ostrowski work on Gauss I without elaboration: instead he offers his own much simpler proof of FTA, the Reich-Smale argument.

The key to parametrizing an algebraic curve is the Implicit Function Theorem. We have need of this theorem in its classical form, but the proof we present is a modern one that is recommended by authorities on (several) complex variable theory.

Our intention is to go through an up-to-date version of Gauss I, driven by the paper  \cite{Gersten-Stallings}. This article is also cited by \cite{MSS} of J. Martin et al., which gives a similar proof using new features, but taking the topological stipulations of \cite{Gersten-Stallings} at face value. The paper \cite{MSS} is directed at combinatorial structures that arise from the interplay of the curves $g(x, y) = 0$, $h(x,y)=0$ when $f(z) = f(x + iy) = g(x,y) + ih(x,y)$.

We redo all the geometry of \cite{Gersten-Stallings} for several reasons. Firstly, it is hardly a good sign to construct regular values of a mapping by using Sard's (non-deterministic) theorem, when the set of critical points (and critical values) is {\it finite} in the first place. One realizes that features of the problem not of interest to the authors are ignored, including certain geometric simplifications. An example is that conceivably $g^{-1}(0) \subset \RR^2_z$ could contain some closed 1-manifold  components, but this is not the case since $g(x,y)$ is a harmonic function satisfying the Maximum Principle. The  ``deep point'' of the authors' proof involves the  ``topology of a 2-cell'' and the Jordan Curve Theorem. We try to make such statements more precise without necessitating a foray into geometric group theory that seems originally adapted to  surfaces of genus greater than $0$.
Sweeping statements are made regarding the convergence of $g(r e^{i \theta})$, ``uniformly with all its derivatives'', but instead of pursuing such results, we find that the first Sections of Ostrowski's paper, the part quoted in Uspensky's book, yield sufficient geometric information to do the job. See \cite{Uspensky}.

The present author admits that by now his (condensed) critique of the \cite{Gersten-Stallings} paper has gone on about as long as that paper itself. One major reason to reconsider this worthwhile article is its use of the ``extended Pre-Image Theorem'', already referred to, that requires a certain mapping to be regular on two different spaces (at the same point). It may be advisable to avoid such an arcane result, especially where we have at hand a visual context of two real variables. You should be able to see the roots emerge as intersections before your eyes! Instead of the ``pre-image as manifold'' result, we use the Implicit Function Theorem as our main tool.

\head The Latest on Gauss's Thesis (1799) \endhead

We consider Gauss's First Proof to fall into the category of proofs based on the open mapping concept. This is because of its critical dependence on the Implicit Function Theorem and the construction of a 1-manifold or plane curve component using an open cover.

Several writings seemed to represent the ``final word'' in describing Gauss's proof, making it sufficiently convincing. This report  will not fully achieve such a goal either. Additionally looking through \cite{Gersten-Stallings 1988}, \cite{Ostrowski 1920}, and \cite{MSS 2002} should provide a good picture of how to carry through Gauss I with modern methods.
Somewhat more difficult than these papers is the original Thesis of Gauss, ``{\it Demonstratio nova Theorematis omnem functionem...}''. A synopsis of essential portions of the  Ostrowski article, which brought the Gauss proof back to a good reputation, is given in the  books of \cite{Uspensky, Appendix I}, or \cite{Fine \& Rosenberger}.

We are going to follow the \cite{Gersten-Stallings} model, but with refinements based on elementary observations about real plane curves. A first aspect is the ``pre-image'' theorem from differential topology, where a given $g^{-1}(y)$ is seen to be a ``one-manifold''.
The desired properties of this point set are then derived from a characterization of an ``abstract'' one-manifold.
We work with concrete curves and arcs, to the extent that ``one-manifold'' becomes superfluous. In particular, we can avoid the full strength of the Pre-Image Theorem. The treatment of this result in well-known books has been debated: we try to step around  issues such as a mapping being transverse to a point set, and also to the boundary of this point set. We are able to avoid use of Sard's theorem with its nondeterministic implications by noting as does Milnor on \cite{TFDV p. 8} that the sets of critical points and  critical values are both {\it finite}.

\newpage
The article \cite{MSS} gives a nice geometric approach different from  \cite{Gersten-Sta\-llings}, by adjusting the component curves instead of the polynomial itself, but takes the needed topological tools for granted.
We do not use the Jordan Curve Theorem at all, falling back on a simpler result that is a prelude to the JCT itself. Also we simplify the final combinatorial step, at the cost of obtaining only one root of the polynomial, not the full contingent of $n$ roots at one time, which the Jordan Curve methods  of \cite{MSS} might achieve.

A. Ostrowski's treatment of a locus of  zeros (a real variety) has been much quoted in its aspect ``toward infinity''. The critical element was where Gauss admitted that he had not proved that an algebraic curve that ``runs into a limited space must run out again''. Ostrowski's clarification of this ``limited'' issue has been completely accepted, but not carried through into any textbook, except by a few drawings \cite{Uspensky, Appendix I}.

Thus it is well to attack the problem from scratch. We have a monic complex polynomial
$$f(z) = z^n + b_{n-1} z^{n-1} + \cdots + b_0,\quad b_j \in \CC$$
which may be rewritten into real and imaginary parts
$$f(x + iy) = g(x, y) + i h(x, y)$$
where $g, h: \RR^2 \to \RR$. We know that $g$ and $h$ are smooth and satisfy the Cauchy-Riemann conditions. One would like to work away from critical points and critical values of $g$ and $h$. Since we will be content to find one zero, with $f(z_0)= 0$, this is not hard to arrange.

Before considering the singularity of the curves $g(x,y) = 0$, $h(x,y)=0$, we use the ``external'' results of Ostrowski, which can be found in greater detail in \cite{Gersten-Stallings}, \cite{MSS}, \cite{Fine \& Rosenberger} and elsewhere. See [Figure D].

\bigskip
{\bf Proposition (Annulus $g$).} {\it There is a real $R > 0$ so  means that the locus of $g^{-1}(0)$ of modulus $r$ over $r \in [R, R+1]$ consists of a quantity $2n$ arcs $\{\gamma_i (t)\}$ where the initial point $\gamma_i (0)$ is $P_i$ and final point is $\gamma_i (1)= P_i'$. Here $|P_i| = R$, $|P_i'| = R+1$},
$$\arg \gamma_i(t) = \frac{(2i+1) \pi}{2n} + \epsilon_i (t), \quad  i = 0, \dotsc, 2n-1.$$

We have that $|\gamma_i (t)|$ is an increasing function, $\epsilon_i(t)$ is smooth with values in $[-1^{\circ}, 1^{\circ}]$ and $\left|\dot{\epsilon}_i (t)\right| < \frac{.01}{R}$. Thus $\arg P_i$ and $\arg P_i'$ are nearly the roots $(2i +1) \pi/2k$ of $\cos n \theta$.\qed

\bigskip
{\bf Proposition (Annulus $h$).} {\it A real value $R > 0$ can be  chosen as above and also so that $h^{-1}(0)$ in the annulus $[R, R+1]$ consists of a quantity $2n$ parametrized arcs $\{\zeta_i(t)\}$ where the initial point $\zeta_i(0)$ is $Q_i$ and whose final point $\zeta_i(1)$ is $Q_i'$ with
$|Q_i| = R$, $|Q_i'| = R+1$. Also $|\zeta_i (t)|$ is an increasing function with
$$\arg \zeta_i (t) = \frac{i \pi}{n} + \delta_i (t),$$
where $\delta_i(t)$ is smooth, deviates from $0^{\circ}$ by at most $1^{\circ}$, and has first derivative $$|\dot{\zeta}_i(t)| \subset \frac{.01}{R} \quad \text{for}\;\; t \in [0,1].\qed$$}

\centerline{
\psfrag{A}{$P_1$}
\psfrag{B}{$Q_1$}
\psfrag{C}{$P_0$}
\psfrag{D}{$Q_0$}
\psfrag{E}{$P_5$}
\psfrag{F}{$Q_5$}
\psfrag{M}{$n=3$}
\psfrag{G}{$P_4$}
\psfrag{H}{$Q_4$}
\psfrag{I}{$P_3$}
\psfrag{J}{$Q_3$}
\psfrag{K}{$P_2$}
\psfrag{L}{$Q_2$}
\psfrag{R}{$R$}
\psfrag{x}{Annulus}
\psfrag{9}{$A_R$}
\psfrag{1}{$P_r'$}
\psfrag{2}{$Q_r'$}
\psfrag{3}{$P_r$}
\psfrag{4}{$Q_r$}
\psfrag{5}{$R+1$}
\psfrag{6}{$P_0'$}
\psfrag{7}{$Q_0'$}
\psfrag{8}{$A_{R+1}$}
\epsfbox{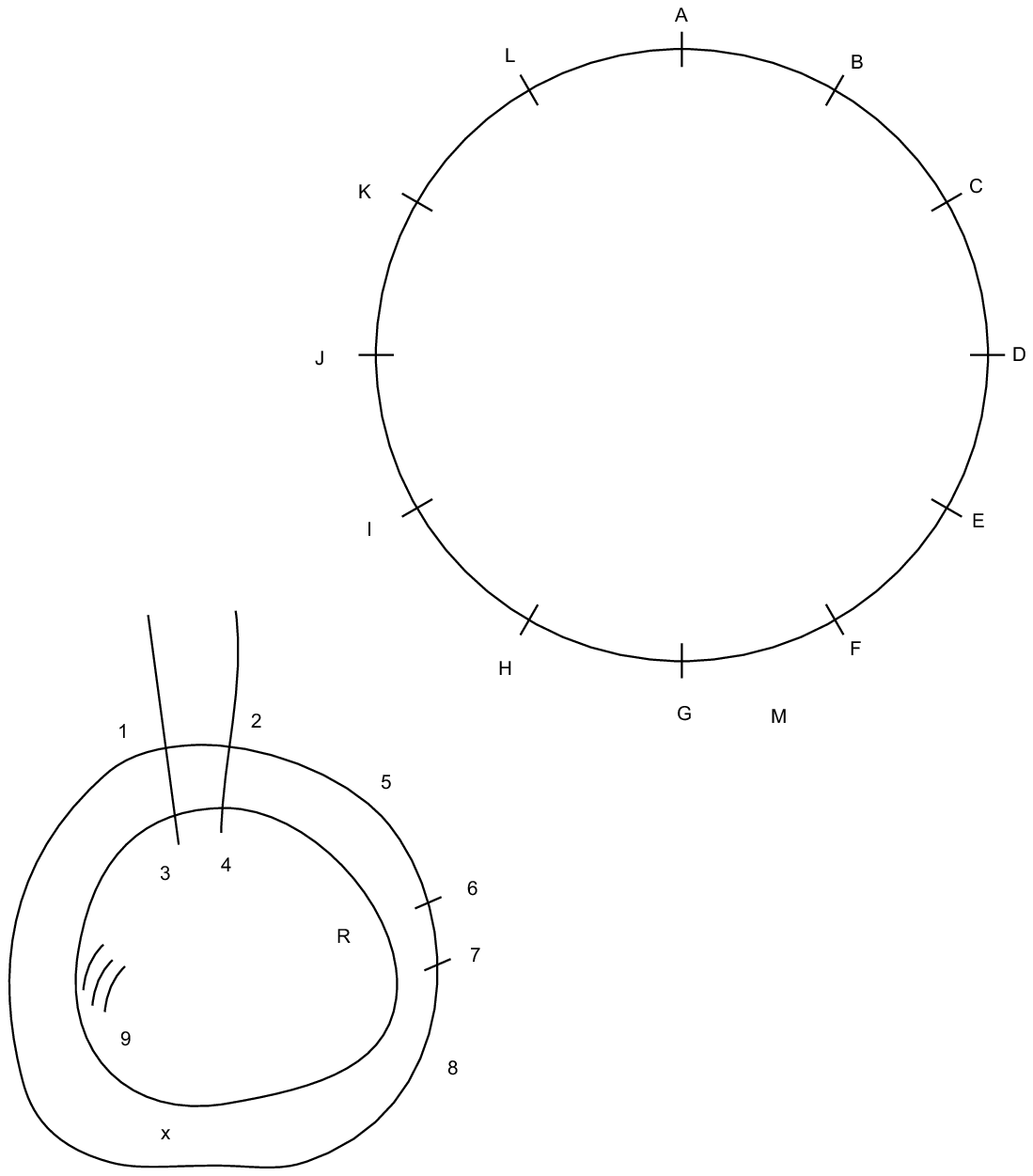}}
$$\text{Figure D}$$

Gauss's conclusion was that the curves defined by $g$ and $h$ meet somewhere within $A_R$, and this point in $\RR^2$ gives a solution to $f(x + i y) = 0$. We prove this using an auxiliary result that is a corollary of the 2-dimensional Brouwer fixed-point theorem. This auxiliary result, ``Crossroads Lemma'', applies to any continuous arc-system (where the arcs may have self-intersections). In fact, the classical Jordan Curve Theorem is a  consequence of the ``Crossroads'' result, see \cite{Maehara}.

\head
W. Walter's Analytic Parametrization
\endhead

With this motivation, we refer back to the Implicit Function Theorem. Local parametrization of an algebraic curve can be dealt with by means of an implicitly defined function such as $F(x, y(x)) = 0$ or $F(x(y), y) = 0$. With $F$ a polynomial, one cannot expect the solution $y(x)$, say to be polynomial. The right category to operate in is that of real analytic functions (power series convergent in some open interval). For example, the ``nodal cubic'' given by $F(x,y) = y^2 -x^2 (x+1)$ has a singularity at $(0,0)$, but can be defined near the Origin by means of two curves (``one-manifolds'') given by $Y(X) = \pm X \sqrt{X+1}$.

At the ``compactifying point''  $(-1,0)$, a separate parametrization, of $X$ in terms of $Y$, should be found in view of $\frac{\partial F}{\partial Y} = 2Y$ which equals  $0$ at $(-1,0)$, even though the curve is smooth here. The square roots in the expressions above can be written as convergent power series. Certainly our starting data, the plane curves that arise as real and imaginary parts of the complex polynomial $P(z)$, form a special case of ``power series'' in two variables.

Thus we use a rather general implicit function theorem, following \cite{Walter 1992}. Let
$$f(x,y) = \sum a_{ij} x^i y^j, \quad f(0,0) = a_{00} = 0, \quad f_y(0,0) = a_{01} \ne 0$$
with $a_{ij}$, $x$ and $y$ belonging to $\RR$.

\bigskip
{\bf Proposition (Implicit Analytic Parametrization).} {\it Suppose that the series defining $f$ converges absolutely for $|x| \leq a$, $|y| \leq b$, with $a, b >0$. Then there  are real numbers $0 < r \leq a$, $0 < s \leq b$ and a power series $w(x)$ converging absolutely for $|x| \leq r$ such that $f(x, w(x)) = 0$ for $|x| \leq r$, $|w(x)| \leq s$ and furthermore $f(x,y) \ne 0$ for all points $(x,y) \in U \times W$, not equal to one of the $(x, w(x))$.}

Here $U \times W$ denotes the rectangular box just constructed.

The uniqueness of solution within the box is critical and may be called ``Walter's Second Uniqueness'', the First being uniqueness merely among analytic solution curves. Actually ``Second Uniqueness'' depends upon carrying through Walter's proof a second time, changing the Banach   algebra of analytic ``germs'' to  a Banach algebra of locally bounded functions. We leave out this additional construction, but instead suggest alternative arguments that are consistent with an ``analytic'' or at least a smooth (differentiable) framework.

\bigskip
{\it  Sketch of proof of Proposition  (see \cite{Walter}).} In the  region of convergence we write $f(x,y) = 0$ in the form
$$Y = \sum_{i,j = 0}^{\infty} b_{ij} x^i y^j := g(x,y)$$
where $b_{00} = b_{0i}  = 0$ and $b_{ij} = -a_{ij}/a_{01}$. We already have the contraction operator that we need. Define $G\,w = g(x, w(x))$ which we will see acts as an operator on a real Banach  algebra $H$. Choose positive real numbers $r,\,s$ according to the recipe
$$\align
B & = \sum_i^{\infty} |b_{i0}| r^i \leq \frac{1}{2}s \\\vspace{0.3cm}
L &= \sum_{i,j}^{\infty} |b_{ij}| r^i j s^{i-1} \leq \frac{1}{2} .
\endalign
$$
Now let $H = H_r$ be the vector space of all functions
$$u(x) = \sum_{r=0}^{\infty} \alpha_k x^k$$
which are absolutely convergent for $x = r$, and define a norm on $H_r$ as
$$\|u(x)\| = \sum_0^{\infty} |\alpha_k| r^k < \infty.$$
It is required to prove that $\|\cdot\|$ on $H_r$ is a legitimate norm, and Cauchy sequences $u_1, u_2, \dotsc, u_f,\dotsc$ of series in $H_r$, converge to a series in $H_r$. Also, with the product $uv$ of the series defining the Banach product, one computes
$$\|uv\| = \sum_k r^k \left|\sum_{i+j=k} \alpha_i \beta_j\right| \leq \sum_k r^k \sum_{i+j=k} |c_{i}| |d_j| = \|u\| \|v\|.$$
Walter gives some basic facts about the Banach algebra $H_r$.
$$\align
\text{\it i)} & \;\|x^k\| = r^k, \;\text{hence}\; \|\bold 1 \| = 1. \\\vspace{0.3cm}
\text{\it ii)} & \;u \in H_r\;\text{implies}\;u^k \in H_r\;\text{with}\;\|u^k\| = \|u\|^k, k = 0, 1, 2, \dotsc \\\vspace{0.3cm}
\text{\it iii)} & \; \text {If } \{u_n\}\;\text{is a sequence in $H_r$ such that}\;\sum\|u_n\| < \infty,\;\text{then}\\\vspace{0.2cm}
& u = \sum u_n \in H_r\;\text{and}\; \|u\| \leq \sum \|u_n\|.\\\vspace{0.3cm}
\text{\it iv)} & \;\text{The integration operator}\; (Iu) (x) = \sum_{k=0}^{\infty} \alpha_k \frac{x^{k+1}}{k+1}\;\text{maps $H_r$ into itself and } \\\vspace{0.2cm}
& \text{satisfies}\;\|Iu\| \leq r \|u\|\;\text{with equality for $u = 1$.}
\endalign
$$
\indent
One may compute that when $\|u\|$, $\|v\| \leq s$ then
$$\|G\, u - G\,v\| \leq L \|u-v\| \leq \frac{1}{2} \|u-v\|.$$
Since $B = \|G (0)\|$, and given $u \in H_r$ with $\|u\| \leq s$,
$$\|G(u)\| \leq \|G(0)\| + \|G (u) - G(0)\| \leq \frac{1}{2} s + L \|u\| \leq s,$$
we see that $G$ maps the closed ball $\|u\| \leq s$ into itself.
By the Banach Fixed-Point Theorem, there must exist a fixed element $\hat{w}$ under $G$, unique for  this property among elements $w \in H_r$.\qed
\bigskip
{\it Proof of Inverse Function Theorem} Let $f: U \to \CC_w$ be analytic and $Df_{|p}$ invertible for $p \in U$. Defining $F : U \times \CC_w \to \CC_w$ by
$$
F(z,w) = f(z) - w \,.
$$
Now $ \dfrac{\partial F}{\partial z} = f'(z)$, so $f'(p) \ne 0$ gives, from the above ``Implicit'' Function Theorem a mapping $g: V \subset \CC_w \to U$ that is locally analytic.  It follows that $F(g(w),w) = 0$ for $w \in V$, in other words $f(g(w)) = w$. But also
$$
F(g(f(z)),f(z)) = f \circ g \circ f(z) - f(z) = 0 \, ,
$$
so $g(f(z)) = z$ for $z \in g(V)$.  Thus we have the two ``inverse'' properties required by the Inverse Function Theorem cited above as a Proposition. \qed

To conclude the Section, we mention Walter's  Second Uniqueness Property, that is, the ``point-wise'' uniqueness of the solution $\hat{w}$ that we found. We repeat the proof above, this time working with the Banach algebra of bounded functions $w: [-r,r] \to \RR$ with norm $\|w\| = \sup  \{|w(x)|: |x| \leq r\}$. This shows as in \cite{Walter} that our (bounded) analytic $\hat{w}$ gives rise to {\it all} the zeros of $u(x,y)$ when $|x| \leq r$, $|y| \leq s$, namely they are exactly the pairs $(x,\hat{w}(x))$. Since we have not covered the proof of Second Uniqueness in detail, those places where it is used in the continuation are given alternate treatment.

\head Regular Values and Curve Singularity \endhead

For the versions of Gauss I carried through on \cite{Gersten-Stallings} and by J. Martin et al. in \cite{MSS}, it is a key point to have both components $g, h$ in $f(z) = g(x,y) + ih (x,y) = 0$ lead to non-singular real algebraic curves $g(x,y)= 0$, $h(x,y) = 0$ valid in a disk $A_R$. Every point $(x,y)$ should be a regular point for both $g$ and $h$, where $(x,y)$ is on the respective curve $g = 0$ or $h = 0$. This avoids self-intersection of any component within $A_R$ of the curve, and for that matter any intersection of two components of $g(x,y) = 0$ (same for $h(x,y)$).

Since for each component $G_0, \dotsc, G_{k-1}, H_0, \dotsc, H_{k-1}$ (as it will turn out), there are only finitely many extrema, we can arrange for the ``coordinate patches'' of this component $\eta_i : [t_b , t_e] \to \RR^2$, $i=0, \dotsc$ to contain at most {\it one} extremum.
There then follows the condition $K$, also referred to as $b)$ in the next Section, which is a key element of the curve construction in \cite{Ostrowski}. At each ``end'' of   $\eta_0$, namely $\eta_0(t_0)$ and $\eta_0(t_1)$ for the endpoints $t_0, t_1$ of the parametrizing interval, the function $\eta_0(t)$ is monotone in both $x$ and $y$ coordinates. Thus definite limits
$$\lim_{t \to t_0^+} \eta_0 (t) = \eta_0^-, \quad \lim_{t \to t_1^-} \eta_0 (t) = \eta_0^+ \;\text{exist.} \tag{K}$$
A similar property holds for all $\eta_j$. One may now use the Implicit Function Theorem to generate $\eta_1: (t_1^*, t_2) \to \RR^2$ on a new interval, centered at $u_1 = \eta_0^+$ and an open set $V_1 \subset \RR^2$ containing $u_1$, where uniqueness of the solution prevails.

The absence of curve singularities is critical to the approach of \cite{MSS} which constructs a beautiful combinatorial structure on the curve components, leading to all   $n$ algebraic roots appearing at once, as intersection points. In our approach we are completely  indifferent to self-intersections and intersections among components. We do need non-singularity (points on the curve are regular for $g$ and for $h$) for one reason: the curve components must have distinct endpoints on the  circle $E_R = \partial A_R$. This will force some component of $G: g(x,y) = 0$ to intersect some component of $H: h(x,y) = 0$, yielding the {\it one} root $z_0 = x_0 + iy_0$ for $f(z)$ that we seek.

Since $g$ and $h$ are harmonic conjugates, the point sets
$$\align
{\Cal S}&= \left\{(a,b) \in \RR^2: g_x (a,b) = g_y(a,b) = 0\right\} \\\vspace{1pt}
{\Cal T} & = \left\{(a,b)\in\RR^2: h_x(a,b) = h_y(a,b) = 0\right\}
\endalign
$$
are the same. In fact this is the ``same'' as $\{z = a+bi\}$ where $f'(z) = 0$, which of course is {\it finite} by elementary algebra, the theory of fields.

We wish both curves to be singularity-free, which means that for any $(a,b) \in {\Cal S} = {\Cal T}$, we have $g(a,b) \ne 0$, $h (a,b) \ne 0$. If there exists $z_0 = a+ ib$ with $f(z_0) = g(a,b) + i h(a,b) = 0$, we have found a root and are done. But it might happen that  $g(a,b) = 0$ and $h(a', b') = 0$ for $a \ne a'$ or $b \ne b'$, $(a', b') \in {\Cal S}$. In that case one or the other of $g$ and $h$ would potentially define a singular curve. Changing $g(x,y) = 0$ to $g(x,y) = \epsilon_1$, by a real constant small in modulus, we may assume that $\tilde{g}(x,y):= g(x,y) -\epsilon_1$ never takes the value $0 \in \RR$ on any $(a,b) \in {\Cal S} =$ set of critical points $\{z_0\}$ of $f(z)$. Similarly we may find $\epsilon_2$ near $0$ such that $\tilde{h} (x,y) := h(x,y) - \epsilon_2$ never satisfies $\tilde{h}(a,b) = 0$ for any $a+ib \, \in$ ``finite singular set of $f(z)$''. Then let $\tilde{f}(z) = f(z) -\epsilon_1 - i \epsilon_2$, which has the same set of critical points as does $f(z)$.

In summary, we wish to modify the  complex equation $f(z) = 0$ so that $g(x,y) = 0$ does not have solutions $(a,b)$ yielding $f'(a + ib) = 0$. The exact same construction applies to $h(x,y)$.

Merely alter $g$ to $\tilde{g}$ by subtracting small  positive or negative $\epsilon_1$. Now the new $\tilde{g}$ might  have acquired a new solution $(a', b')$ where $f'(a' + ib') = 0$. In that case push all $\tilde{g}$ to $\hat{g}$ by adding to $\tilde{g}$ a real constant $\epsilon_1'$, smaller in modulus than $\epsilon_1$, so by now we have avoided both ``critical'' solutions $(a,b)$ and $(a',b')$. After finitely many steps we have (re-using notation) $\hat{g}(x,y) = g(x,y) + \hat{\epsilon}$ where $\hat{g}(x,y) = 0$ contains no singularities. Again by closure of $f(z)$, given that $f(z) = 0$ also has {\it no} solution, we construct $\epsilon, \delta$ where $f(z) = \epsilon + i \delta$ has {\it no} solution, and its constituent real harmonic curves $g= 0$, $h=0$ have only regular points.

Admittedly the somewhat lengthy argument above is covered by  \cite{Gersten-Sta\-llings} in one sentence. But the authors did not make explicit the need to assume, for the purposes of their argument, that both harmonic curves are non-singular.

We just established that there is a sequence $\{\epsilon_1^k\}$ converging monotonically in modulus to  $0 \in \RR$ (where $k$ is an {\it index}), such that $g(a,b) = \epsilon_1^k$ {\it never} has a  solution in ${\Cal S}$. Also we have a sequence $\{\epsilon_2^k\}$ converging monotonically in modulus to $0 \in \RR$ such that $h(a,b) = \epsilon_2^k$ never has a solution $(a,b) \in {\Cal S} = {\Cal T}$ either. We claim that if $f(z) = 0$ has no solution at all, then neither does $f(z) = \epsilon_1 + i \epsilon_2$, for values $\epsilon_1, \epsilon_2$ arbitrarily close in to $0 \in \RR$. If such a convergent sequence did exists, with solutions $z_k$
$$f(z_k) = \epsilon_1^k + i \epsilon_2^k,$$
the solutions would be bounded and a convergent sub-sequence of $\{z_k\}$ would lead to $f(z_k) = 0$. Thus the Reich Principle shows that we can reduce the problem of existence of a root for $f$ to one where the two real curves $g(x,y) = 0$ and $h(x,y) = 0$ have no singularities in $\RR^2$.

\newpage

With these choices we now have in $A_{R+1}$, that $g^{-1}(0)$ is a ``smooth 1-manifold'', consisting of several arcs with no intersections, and $h^{-1}(0)$ is also a ``smooth 1-manifold'' composed of  non-intersecting arcs.

Furthermore $g^{-1}(0) \cap h^{-1}(0)$ is empty unless some $x + iy$ in the intersection solves $f(x + iy) = 0$. We write $g^{-1}(0)$ instead of $g^{-1}(\epsilon)$ as in \cite{Gersten-Stallings} as we take it that the ``$\epsilon$ modification'' to the original polynomial function has already been carried through. The following Section will use the Implicit Function Theorem to describe the arc structure of $g^{-1}(0)$ {\it within} the disc $A_R$. The boundary points of $g^{-1}(0)$ and $h^{-1}(0)$  on $E_R$ or $E_{R+1}$ lie on a combinatorial configuration that eventually will contradict $g^{-1}(0) \cap h^{-1}(0) = \emptyset$, and we  will produce a solution to $f(\hat{z}) = 0 \in \CC_w$.

\head Inside  the Disk $A_R$ \endhead

We recall the construction of points $\{ P_i,Q_i\}$ on $\partial A_R$ that constitute ``inner end points" of the arcs $\gamma_i$ and $\zeta_i$, referring again to [Figure D]. Now the restricted sets $\gamma = A_R \cap g^{-1}(0)$ and $\zeta = A_R \cap h^{-1}(0)$ are defined as plane curves, and we wish to characterize those connected arcs in $A_R$ that represent a continuation at $P_i$ or $Q_i$ of a given ``exterior" arc $\gamma_i$ or $\zeta_i$.
Let us concentrate on the case $\gamma$, the case of $\zeta$ will be similar.  We have ``endpoints" $\{P_i\}, \, i=0, \dotsc , 2n-1\,$. Other points $E \in A_R$ of interest are those where $x(\gamma)$ or $y(\gamma)$, attains a local maximum or minimum. by B\'ezout's Theorem, the cardinality of these {\it extremal} points is certainly no greater than $n(n-1)$, it is {\it finite}.

Let $P$ be one of the $\{ P_i\}$. By Implicit Parametrization above (Walter's theorem), we may find a one-sided analytical arc $\eta : [0,1) \to A_R$ expressing $\eta(0) = P$, $\eta(t) \in \text{Int}A_R$ for $t > 0$.  Hence $\eta$ can be considered as a diffeomorphism from the half-open interval to an arc $\eta_P$.  This arc can be chosen not to intersect any of the extremal points $E$. The same procedure is followed at every boundary point from $\{P_i\}$, and also at all the extremal points from $\{E\}$, save that in the latter case we end up with a two-sided open arc $\eta_E$ whose image has $E$ in its interior.  What one must now do is to extend these arcs to obtain the curve-components $\gamma_{P_i},\, \zeta_{Q_i}$ that connect pairs $P_i$ and $P_{\pi_i}$, $Q_j$ and $P_{\pi_j}$ on $\partial A_R$.

For each arc, Property K above applies and the arc may be extended from one limiting end point or the other, or both (in the case of an extremal location of the type E), until its closure contains an extremum or boundary point. In practice, we extend (by the Implicit Function Theorem above) only until an overlap {\it of arcs} occurs.  Thus the extending patch runs into another patch that originated from some extremum or boundary point from $\partial A_R$.  Continuing exhaustively in this manner, there results a collectio of patches $\eta_j : (0,1) \to A_R$, given by $\eta_j(t) = [x(t),y(t)]$, where locally the functions $x(\cdot)$ and $y(\cdot)$ are given as convergent power series in $t$. See [Figure E]. 


\bigskip
\centerline{
\psfrag{1}{$t_0$}
\psfrag{2}{$t_1^*$}
\psfrag{3}{$t_1$}
\psfrag{4}{$t_2$}
\psfrag{5}{$\sigma(t_0)$}
\psfrag{6}{$\sigma (t_1)$}
\psfrag{7}{$\sigma(t_1^*)$}
\psfrag{8}{$\sigma (t_2)$}
\psfrag{R}{$\RR$}
\epsfbox{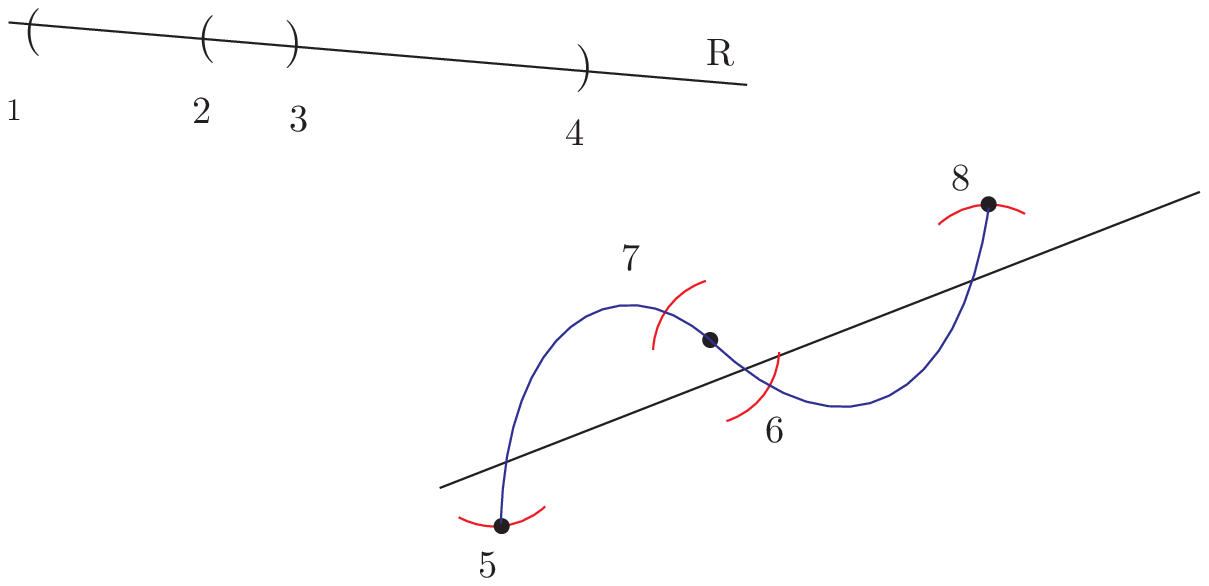}}
$$\text{Figure E}$$

Three criteria for $\{\eta_j\}$ hold:
\roster
\item"(a)"
  $\eta_j$ has at most one $y$-extremum (maximum or minimum) in the parametrization $[x, y(x)]$ and at most one $x$-extremum in the parametrization $[x(y), y]$.

\item"(b)" Unique limits exist for $x(t)$, $y(t)$ as $t \to t_0^-$ and  $t \to t_1^+$.

\item"(c)"  The $\{\eta_k\}$ are ordered linearly with a non-empty overlap  $\text{Im} (\eta_i) \cap \text{Im} (\eta_j)$ only when $i = j+1$ or $j = i+1$.
\endroster

Property (a) follows from the construction, which never allows a new extremal point into the patch that is being extended from an originating extremum, or boundary point. This ensures that a limit at either end $\eta(0)$ or $\eta(1)$, is guaranteed to exist.  Property (b) is the same as Property K mentioned earlier, and results from $x(\eta),\,y(\eta)$ being {\it monotone} functions of $t$ near the endpoint os $I=[0,1]$, whether $\eta$ is defined on an open interval or a semi-closed interval. Property (c) is the subject of the remainder of the Section.

{\bf Discussion of $\eta$ conditions.} Since our problem relates to the topology of curves in the plane, the intersection of $g^{-1}(0)$ and $h^{-1}(0)$, we may adjust the coordinate system to gain any advantage through Algebra. In particular we want $g(x,y)$ as a polynomial form to contain no single-variable factors $b(x)$ or $c(y)$: these would present curve components parallel to an axis. This being given, the number of extrema on $g^{-1}(0)$ should not be greater than $2n (n-1)$, as follows from B\'ezout's Theorem \cite{G. Fischer, Section 3.2}.

Since $g$ and $h$ are continuous functions, we have that $g^{-1}(0)$ and $h^{-1}(0)$ are closed subsets of $A_R$. The connected component of $g^{-1}(0)$ containing $P_0$ is constructed as above by a sequence of arcs $\{\eta_i\}$ coming  from Walter's Implicit Function Theorem. The  arcs will eventually exhaust the allowable finite number of extrema. The ``final'' arc $\eta_w$ will either ``stop suddenly'' in the interior of $A_R$, or meet $\partial A_R$. By ``final arc'' we may mean a convergent sequence of monotone arcs. In either case one can construct a global parametrization of  that part of the component $G_0$ reached to this point, as a concatenation, leaving in mind overlap of  the local parametrizations $\{\eta_i\}$ coming from Implicit Function Theorem. In the case where convergent ``ends'' of a sequence $\eta_k, \eta_{k+1},\dotsc$ converge to  $u^* = (x^*, y^*) \in \RR^2$, we may take $u^*$ as the center of a new local parametrization $\eta^* : (t_{\infty}, t_{\infty} + \epsilon) \to \RR^2$.  See [Figure F].

The other possibility is that a $\eta_j$ intersects the image of a previous $\eta_i$, $i < j$ or, a sequence  $\eta_j, \eta_{j+1},\dotsc$ comes arbitrarily close to an image point of $\eta_i, i < j$. Specifically, the open set $V_i$, ``domain of uniqueness'' can be intruded on by patches that were generated subsequent to $\eta_i$.

\bigskip
\centerline{\psfrag{1}{$P_0$}
\psfrag{2}{$\sigma_0$}
\psfrag{3}{$G_0$}
\psfrag{4}{convergent sequence}
\psfrag{5}{new arc}
\psfrag{e}{$\eta^*$}
\psfrag{6}{$u^*$}\epsfbox{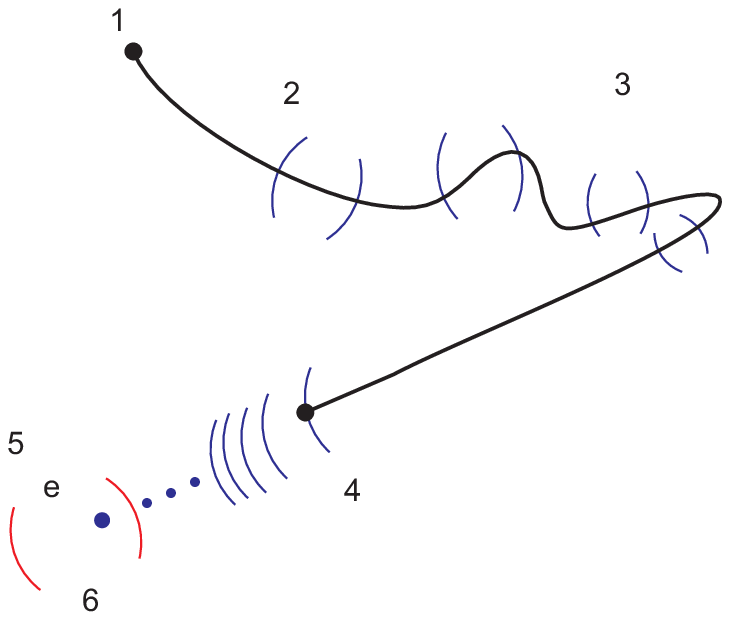}}
$$\text{Figure F}$$

\bigskip
Walter's Second Uniqueness result, part of the analytic Implicit Function Theorem, rules out  such behavior. If $D_i = \{x,y : r_1 < x < r_2,\, s_1 < y < s_2\}$ is the domain of uniqueness for the patch $\eta_i$, then the  only values $(x,y)$ in $D_i$ that satisfy $g(x,y) = 0$ are the values $\eta_i(t) = \left(x(t), y(t)\right)$ for $t$ in the parametrizing interval $\left(t_i^*, t_{i+1}\right)$. See [Figure G].

As previously remarked, this part of Walter's Theorem requires consideration of a Banach algebra larger than ``locally convergent power series'', namely ``locally bounded functions''. It would be good  to prove this uniqueness (a double point or crossing is an algebraic singular point) without leaving
the category of power series. For example, if $\eta_j$  were to merge with $\eta_k$ with an infinite order of tangency, all higher order derivatives at $u$, namely $\dfrac{dy}{dx},\, \dfrac{d^2y}{dx^2},\dotsc$ are equal for the two curves. Thus by uniqueness of analytic solution the curves are equal in a neighborhood of $u$. But $u$ was assumed to be the first point for the  parametrization that the curves meet (the curves are  topologically closed) which gives a contradiction.

\bigskip
\centerline{
\psfrag{1}{$\eta_k$}
\psfrag{2}{$\eta_{k+1}$}
\psfrag{3}{convergent}
\psfrag{4}{$\eta_0$}
\psfrag{5}{$u^*$}
\psfrag{6}{$\partial A_R$}
\psfrag{7}{$P_k$}
\psfrag{8}{$\eta_w$}
\psfrag{9}{$A_R$}
\epsfbox{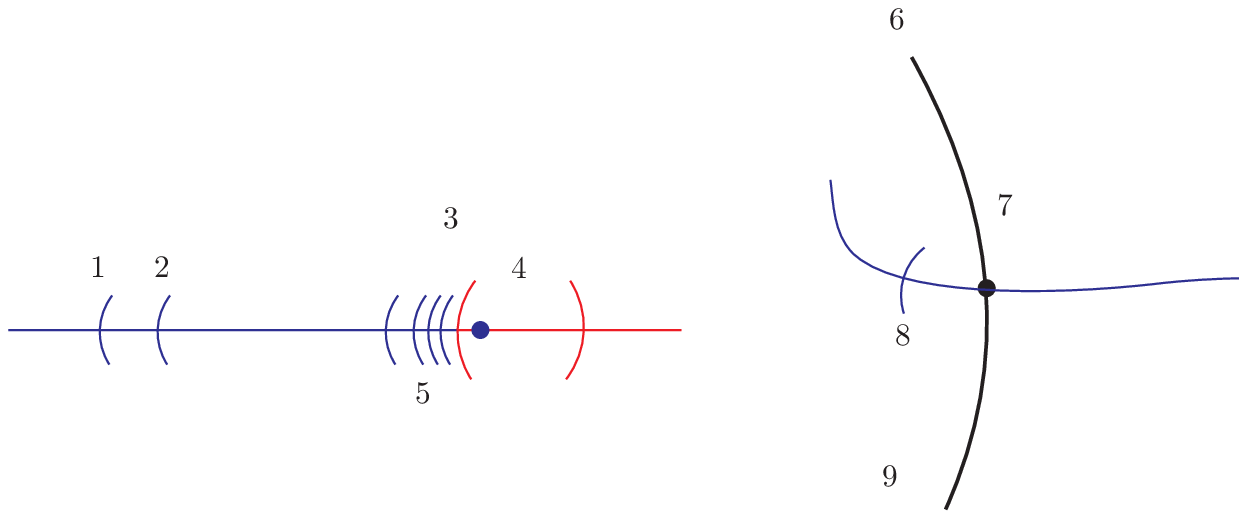}}
$$\text{Figure G}$$

\bigskip
\bigskip
\centerline{
\psfrag{u}{$u$}
\psfrag{1}{$\eta_j$}
\psfrag{2}{$\eta_k$}
\epsfbox{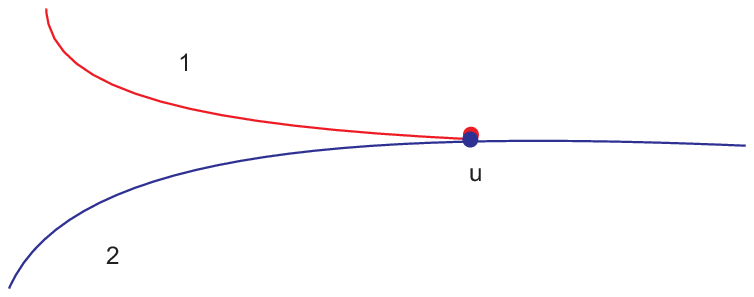}}
$$\text{Figure H}$$

\bigskip

If on the other hand, $\eta_j$ and $\eta_k$ differ at $u$ in some power of tangency, there are formulas that specify this ``slope'' or tangency, and there is  no  leeway for solutions to the relation $g(x,y) = 0$ locally. For example if $y = \sigma(x)$ is the solution at regular point $x, y$ the slope there is given by the well-known                        formula
$$\left.\dfrac{dy}{dx}\right|_{p = (x,y)} = \dfrac{d \sigma (x)}{dx} = \dfrac{-g_x(x,y)}{\,\,g_y(x,y)}$$
\indent The formula for the second derivative is
$$\dfrac{d^2y}{dx^2} = \dfrac{\left(-g_y^2 g_{xx} + 2g_x g_y g_{xy} - g_x^2 g_{yy}\right)}{g^3_y}$$
where
$$g_{xx} = \dfrac{\partial^2 g}{\partial x^2}, \; g_{xy}=\dfrac{\partial^2 g}{\partial x dy}\;\;\text{etc.}$$
There are formulas for all order derivatives, valid as long as $g_y \ne 0$. This shows the Taylor ``jet'' or ``germ'' at $P$ is completely determined by $g(x,y)$ as long as $g$ is regular (surjective) at $P$.  See [Figure H].

The above considerations have an essential consequence. Though we noted that it is not vital for the rest of the proof whether $G_0$ has any ``self-intersections'' or whether $G_i$ intersects $G_j$ for $i \ne j$, it is essential that the starting node $P_{i_{b}}$ of $G_i$ be distinct from the ending node $P_{i_{e}}$, and that this pair be {\it disjoint} from any pair $P_{j_b}, P_{j_e}$ for $i \ne j$.
We essentially did show that no self-intersection, mutual crossings or mergings between $G_i, G_j$ can occur, which  is key to the ``basketball'' argument in \cite{MSS}.

A possible drawback of the reasoning about arcs given above is that either one must work through a different ``Walter'' Uniqueness  argument in a new category (bounded functions) or one must apply background knowledge about germs and jets of convergent power series. An alternative will now be sketched, that keeps us in the smooth category which is familiar to many. Taking by the argument about ``finitely many extrema'' of all the component curves $G_0, G_1, \dotsc, G_{k-1}$ (at least we suspect that they are curves) we may look at an intersect or ``merger'' point isolated in a rectangular box [Figure I]:

\centerline{
\psfrag{1}{$\Gamma$}
\psfrag{2}{$\Delta$}
\psfrag{3}{box $B$}
\psfrag{g}{$g$}
\psfrag{R}{$\RR$}
\psfrag{w}{$w$}
\psfrag{P}{$p$}
\epsfbox{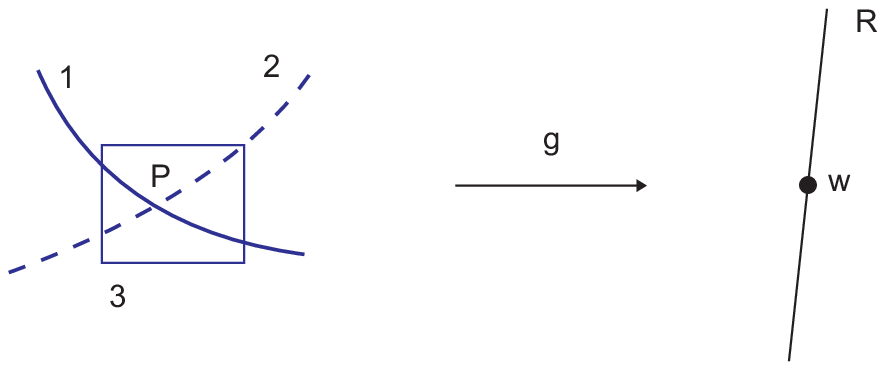}}
$$\text{Figure I}$$

\bigskip

Now since $p$ is regular for $g(x,y)$, we can apply the Local Submersion Theorem of differential
topology \cite{Guillemin \& Pollack, Section 1.4}. This is  proved directly from an Inverse Function Theorem that is available to us. Local submersion tells us that there is a diffeomorphism $\psi: B \to T$ where $T$ is another box but $\psi(\Gamma) = L$, where $L$ is a horizontal segment. See [Figure J].

\bigskip
\centerline{
\psfrag{G}{$\Gamma$}
\psfrag{D}{$\Delta$}
\psfrag{B}{$B$}
\psfrag{R}{$\RR^2$}
\psfrag{P}{$p$}
\psfrag{T}{$T$}
\psfrag{L}{$L$}
\psfrag{a}{$a$}
\psfrag{b}{$b$}
\psfrag{f}{$\psi$}
\epsfbox{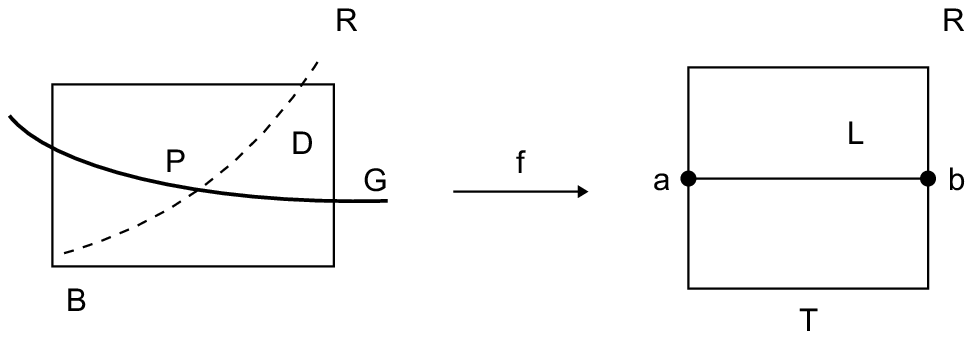}}
$$\text{Figure J}$$

\newpage
But applying the same theorem to $\Gamma + S$, $S$ a subset of the image of $\Delta$ (``the other arc'') gives another diffeomorphism $\varphi: B \to T$ where $\varphi (\Gamma) \subseteq L$, but also $\varphi (g) \in L$  where $g \in S$. The diffeomorphism $\lambda = \psi \varphi^{-1}: T \to T$ takes $L \to L + S$. The map $\lambda$ is monotone on $L$ and maps $\psi(p)$ onto $\varphi(p)$. Points in $S$ converge to $p$; therefore some $\lambda^{-1}(s)$ lies on the interval $\left[\lambda(a), \lambda (b)\right]$. But then $\lambda^{-1}(s)$ is not attained on $t \in [a,b]$, contradicting the Intermediate Value Theorem.

Again, the ``terminal'' boundary point of $K = \text{Im}\,\sigma$ must be some $P_m$ distinct from $P_0$ for the reasons just propounded. That is, $\sigma_w$ would have to merge with $\sigma_0$ at a previous coordinate $\sigma_i(t')$, or meet $P_0$ directly from inside $A_R$. The $\sigma_w(t)$ values near $P_0$ would provide ``extra solutions'' to $g(x,y) = 0$ that are ruled out by Walter's Second Uniqueness Theorem. Alternatively one can show the same, that $G_i$ has two distinct endpoints of $\partial A_R$, and these are distinct from those of all other $G_j$, by means of the derivative formulas and analytic uniqueness, or by the Local Submersion Theorem of differential topology.

We recapitulate the situation regarding  algebraic arcs inside a closed disc $A_R$. We quote C.F. Gauss (see \cite{Smale}), ``an algebraic curve can neither suddenly be interrupted... nor lose itself after an infinite number of terms''.

From our point of view, the curve cannot ``suddenly be interrupted'' unless       $\partial A_R$ is reached, since an extension  of the growing arc can always be found at any limit point such a ``$u$'' discussed above. The curve cannot ``lose itself'' into oblivion like a logarithmic spiral, since the number of $x$- or $y$- extrema would have no bound. Arguments from compactness were not available until after Gauss's time, but such a proof using B\'ezout's theorem would have been at hand.

So, according to Gauss, there remains the possibility that the curve ``runs into itself'',
which we could rule out since we have enforced non-singularity of the curve components. There remains only ``runs out to infinity in both directions'' (at distinct angles), which means that each topological component such as $G_i$ has two boundary points on $\partial A_R$.

The admission by Smale, Master of the high-dimensional Universe, that ``it is a subtle point even today'' why a real algebraic component $g^{-1}(0)$ cannot enter $A_R$ without leaving, makes one wonder whether all similar issues have been cleared up for ``3-folds in projective $N$-space'' and so forth.

Pulling together the various pieces, we apply the process given above to all components of $G = g^{-1} (0) \cap A_R$ and components of $H = h^{-1}(0) \cap A_R$. We find, as is discussed in \cite{Gersten-Stallings}, \cite{Uspensky} and \cite{MSS} that there are $n$ arcs $G_{0}, \dotsc, G_{n-1}$, parametrized by $\{\sigma_i\}$, and $n$ arcs $H_0, \dotsc, H_{n-1}$, parametrized by $\{\tau_j\}$, connecting up the $P_0, \dotsc, P_{2n-1}$ and $Q_0, \dotsc, Q_{2n-1}$ respectively. In [Figure K] we see the ``matching'' partially defined by $[0] \leftrightarrow [k]$ for $P$ and $[0] \leftrightarrow [m]$ for $Q$.

In the previous Section we saw that the collection of $P$-arcs $\{\sigma_i\}$ in $A_R$, each corresponding to a component $G_i$, were disjoint by the smoothness of the overall algebraic curve $g(x,y) = 0$. Similarly the arcs $\{\tau_j\}$ corresponding to the components $\{H_i\}$,  whose endpoints are $\{Q_f\}$ do not intersect. The goal now is to show that some arc $\sigma_j$ must meet some arc $\tau_k$ within $A_R$.

\newpage
\centerline{
\psfrag{A}{$A_{R+1}$}
\psfrag{B}{$Q_m'$}
\psfrag{C}{$P_k$}
\psfrag{D}{$Q_m$}
\psfrag{E}{$\sigma_0$}
\psfrag{F}{$\tau_0$}
\psfrag{G}{$P_k$}
\psfrag{H}{$P_k'$}
\psfrag{I}{$P_0$}
\psfrag{J}{$P_0'$}
\psfrag{K}{$\gamma_0$}
\psfrag{L}{$Q_0$}
\psfrag{M}{$\xi_0$}
\psfrag{N}{$Q_0'$}
\epsfbox{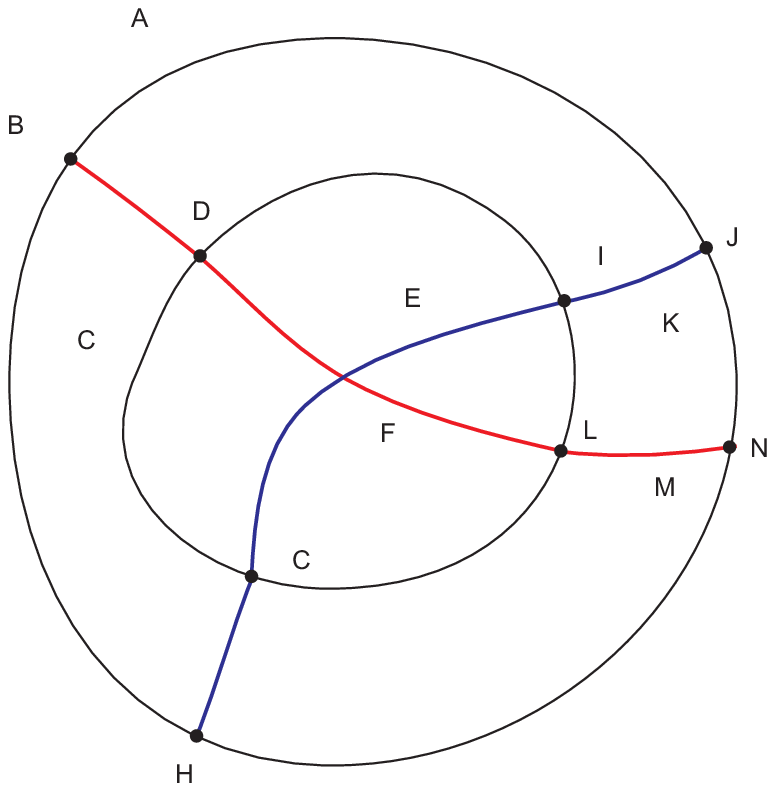}}
$$\text{Figure K}$$

\bigskip
We give the topological part of the short remaining argument.

\bigskip
{\bf Proposition [Maehara].} {\it Suppose that a continuous arc $\sigma$ on $A_R$ has distinct endpoints (nodes) $E, F$ separated by nodes $A,B$, where $E = [2 e]$, $F = [2 f]$, $A = [2a+1]$, $B = [2b+1]$ are the distinct boundary points of $\tau$. Then $\sigma$ and $\tau$ have a common point (non-empty intersection) within $A_R$.} \qed

\bigskip
\centerline{
\psfrag{1}{$A_R$}
\psfrag{2}{$B = Q_b$}
\psfrag{3}{$F = P_f$}
\psfrag{4}{$A = Q_a$}
\psfrag{5}{$E = P_e$}
\epsfbox{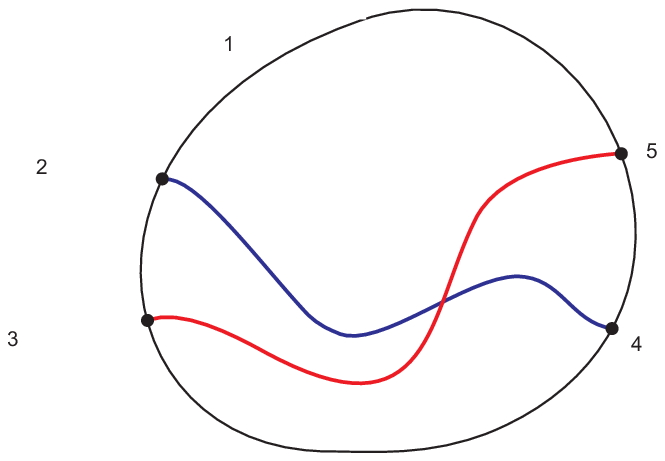}}
$$\text{Figure L: Maehara Crossroads Theorem}$$

\bigskip
{\bf Remarks.} Note that $\sigma$ and $\tau$ in this statement are not required to be {\it simple} or smooth arcs, but each boundary $\partial \sigma$, $\partial \tau$ lies on $\partial A_R$ and consists of  two points. ``Separated'' means that on the circle, reading counterclockwise the indicated nodes similar to the following, see \cite{MSS}.
$$EQP\;Q\quad PQ \;\ldots\; APQP\;Q \;\ldots \; FQPQPQ\ldots BPQP\ldots$$
or
$$EQPQ \;\ldots\; BPQP\;\ldots\; FQPQ\;\ldots\; APQP\;\ldots,$$
or a {\it cyclic permutation} of same.

On the other hand the configuration $EPQPFAQPQB$ does not satisfy the hypothesis. In this case it is possible to choose $\sigma$ and $\tau$ that do not intersect. See [Figures L \& M].

\bigskip
\centerline{
\psfrag{A}{$A$}
\psfrag{B}{$B$}
\psfrag{F}{$F$}
\psfrag{E}{$E$}
\psfrag{t}{$\tau$}
\psfrag{s}{$\sigma$}
\epsfbox{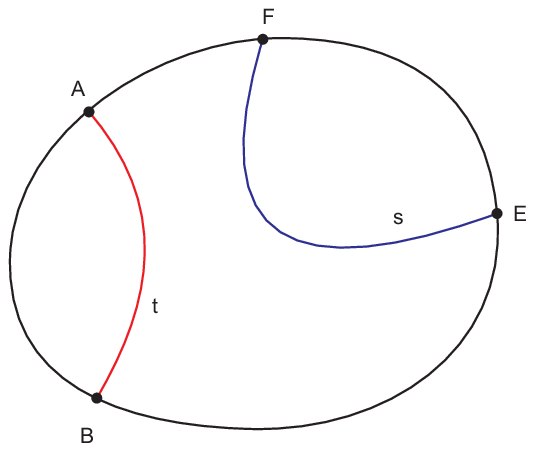}}
$$\text{Figure M}$$

\bigskip
A generalization of Maehara's result, deriving from the Theorem of Poincar\'e-Miranda, is discussed in the Appendix. The article by three authors on ``basketball configurations'' \cite{MSS} shows more strongly that {\it every} $\sigma_i$ is matched to  {\it one} $\tau_{j_i}$ with which it has {\it exactly one} intersection point. Their proof uses non-self-intersection of $G$ and $H$ and a less elementary topological fact, the Jordan Curve Theorem (in its form applying to smooth curves).

\head Sector Matching by Harmonic Components \endhead

We review notation that has already been used, and is consistent with the treatment in \cite{Uspensky, Appendix I}, and similar to that of \cite{MSS}. Consider $2 N$ non-negative integers in an ordered interval
$$\{2 N\} = [0, 1, \dotsc, 2 N-1].$$
Choosing say $N=5$ and $A = 3$, $B = 7$ we obtain two {\it sectors} of $\{2N\}-(A,B)$, namely $S_1 = [4, 5, 6]$ and $S_2 = [8, 9, 0, 1, 2]$ where we see that the  integers were actually in sequence $\mod (10)$. The model to keep in mind is the circle $E_R = \partial A_R$ described earlier where $P_{0}, \dotsc, P_{2N-1}$ were given                in counter-clockwise order. We may say that 4 and 6 are in the {\it same sector} $S_1$ for  $\{2N\}-(A,B)$ but that the ``matched pair'' $(A, B)$  {\it separates} nodes 5 and 0. See [Figure N].

Getting to the case of interest, let $N = 2n$ where  $n$ is the degree of our original polynomial $P(z)$ (or $f(z)$). In this case our geometric labels look like
$$Q_0 \sim [0] \quad Q_1 \sim [2] \cdots Q_{N-1} \sim [4n-2]$$
\vskip-0.6cm
$$P_0 \sim [1] \quad P_1 \sim [3] \cdots P_{N-1} \sim [4n-1]$$
in the above description.

\newpage
The pairing of the boundary by arcs $\sigma_j$ corresponding to $P_j$, and the boundary points of arcs $\tau_k$ corresponding to node $Q_k$ gives a {\it matching} (fixed-point free involution) of the  $\{P_j\}$, and of the $\{Q_k\}$ respectively. Given some arc $\sigma$, its distinct end-nodes form two sectors $I$ and $II$
$$A \sim [k], \; [k+1], \dotsc, \; B \sim [m].$$
Hence, $A$ and $B$ are ``$P$-nodes''. Suppose that the number of $Q$-nodes in sector $I$ is {\it odd}. See [Figure O].

\centerline{
\psfrag{1}{Sector $I$}
\psfrag{2}{Sector $II$}
\psfrag{3}{$P_k\;Q_k\;P_{k+1}\;Q_{k+1} \cdots P_{k+5}$}
\psfrag{4}{$Q_{k+5} \cdots$   }
\epsfbox{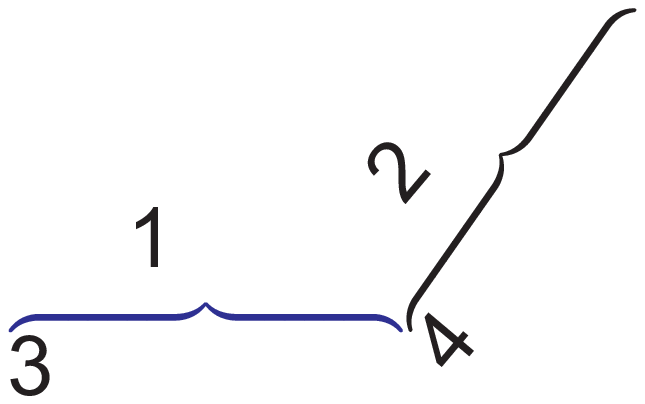}}
$$\text{Figure N}$$

\vfill
\centerline{
\psfrag{S}{Sector $I$}
\psfrag{T}{Sector $II$}
\psfrag{A}{$P_2$}
\psfrag{B}{$Q_1$}
\psfrag{C}{$P_1$}
\psfrag{D}{$Q_0$}
\psfrag{E}{$P_0 \sim [0]$}
\psfrag{F}{$Q_7 \sim [15]$}
\psfrag{G}{$P_7 \sim [14]$}
\psfrag{H}{$Q_6$}
\psfrag{I}{$P_6$}
\psfrag{J}{$Q_5$}
\psfrag{K}{$P_5$}
\psfrag{L}{$Q_4$}
\psfrag{M}{$P_4$}
\psfrag{N}{$Q_3$}
\psfrag{O}{$P_3$}
\psfrag{P}{$Q_2$}
\psfrag{Q}{$\tau_1$}
\psfrag{R}{$\sigma_0$}
\epsfbox{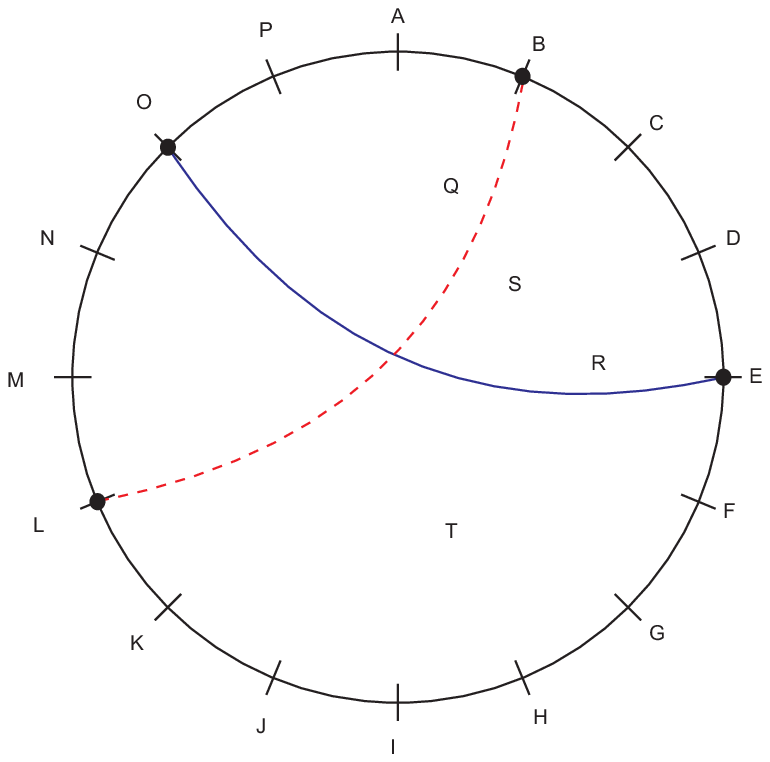}}
$$\text{Figure O: $(P_0, P_3)$ gives
$\text{Sector}\;I = [Q_0, P_1, Q_1, P_2, Q_2]$ and}$$
 $$\text{Sector}\;II = [Q_3, P_4, Q_4, P_5, Q_5, P_6, Q_6, P_7, Q_7]$$
hence separates $Q_1$ from $Q_4$. If $\sigma_0$ has boundary $P_0, P_3$,
then $\tau_1$ with boundary $Q_1, Q_4$ should intersect $\sigma_0$.

\newpage

  Then one of the nodes $Q_f$ in Sector $I$ must match (be attached by a $\tau$ arc) a node $Q_f'$ of Sector $II$. By the Proposition of Maehara, this $\tau$-arc must intersect the original $\sigma$-arc in the interior of $A_R$. So we would have a common solution $(a,b) \in \RR^2$ for  $g(a,b) = h(a,b) = 0$.

  On the other hand, if the number of $Q$ nodes in Sector $I$ is even, we may assume that all of their arc-pairings occur {\it within} Sector $I$, else we have a $\tau'$ that must meet $\sigma$ as before. Any such $Q$-pairing, call it $\tau''$, forms new sectors labeled $III$ and $IV$, one of which, say $III$, lies entirely within the $P$-sector $I$, hence is  strictly smaller in cardinality. Now we are interested in the $P$-nodes of Section $III$. As always under this construction (when Sector $\Lambda'$ ends up strictly contained in Sector $\Lambda$), there is at least {\it one} $P$-node in Sector $III$. If the count of these $P$-nodes is odd, a pairing arc $\sigma'$ must arise that meets $\tau''$ in a solution point. If we are still not finished, interchanging the r\^oles of $P$ and $Q$, $g$ and $h$ and so forth, leads by induction to a basic case of a singleton $P$- or $Q$-node that must be paired outside its sector, leading to a solution point.

\head Appendix: The ``Crossroads Theorem'' in Higher Dimensions \endhead

Maehara's `Crossroads' result generalizes from arcs in a disc or square, to ``hypercurves'' of complementary dimension, transverse in a cube of the ambient dimension.

If $I = [-1,1]$ is the closed double interval, we define $I_k \subseteq I^n$ by
$$I_k = \left\{|x_1| \leq 1, \dotsc, |x_k| \leq 1, \; x_{k+1} = 0, \dotsc, x_n = 0\right\}$$
and similarly $\hat{I}_l \subseteq I^n$ by
$$\hat{I}_l = \left\{x_1 = 0,\dotsc, x_{l-1} = 0, \;|x_l| \leq 1, \dotsc, |x_k| \leq 1\right\}.$$
Suppose we have mappings $g: I_k \to I^n$, $h: \hat{I}_{k+1} \to I^n$ satisfying
$$g_{|\, \partial I_k} = id_{\partial I_k}, \qquad
h_{|\, \partial \hat{I}_{k+1}} = id_{\partial \hat{I}_{k+1}}.$$
Then we have

\bigskip
{\bf Proposition (Generalized Crossing Theorem).} {\it In this case there exist $s \in I_k, t \in \hat{I}_{k+1}$ such that
$$g(s) = h(t).$$
In other words, the image $G = g(I_k)$ meets the image $H = h (\hat{I}_{k+1})$ in at least one point.}

\bigskip
Maehara's result is when  $k=1$, $n=2$. For example, consider an arc (a `path') in $I^3$ from $A$ to $Q$ (in red), and a surface within $I^3$ whose boundary is the ``equator'' in blue, $BCDE$. Then the path and the surface must meet within {\it closed} $I^3$. See [Figure P].

The proof is an immediate application of Miranda's Theorem, a version of the Brouwer Fixed-Point Theorem originally proposed by Poincar\'e. See \cite{Miranda}, \cite{Vrahatis}. \qed

We regard Brouwer's FPT as a tool to be employed without hesitation. Proofs of the equivalent ``non-retraction theorem'', due to Y. Kannai (see \cite {Flanders}), C.A. Rogers, and Milnor-Asimov are elementary and lucid.
Any of these approaches leads to a modern proof of the Poincar\'e-Miranda Theorem, the generalized Crossroads Theorem (``topological transversality'') and a new proof of the result of Maehara, which he uses in turn to obtain a short proof of the Jordan Curve Theorem.

\bigskip
\centerline{
\psfrag{E}{$E$}
\psfrag{Q}{$Q$}
\psfrag{D}{$D$}
\psfrag{B}{$B$}
\psfrag{C}{$C$}
\psfrag{R}{$A$}
\epsfbox{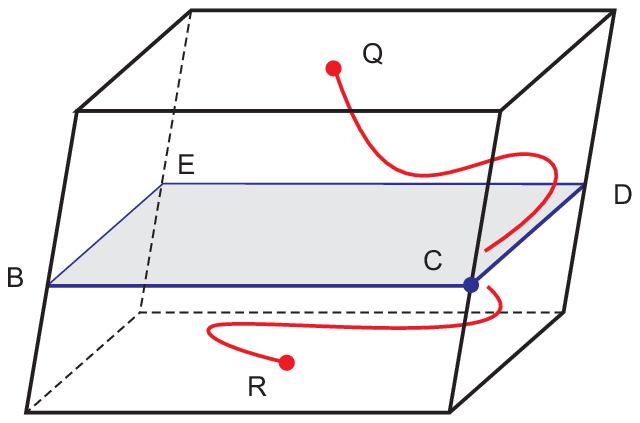}}
$$\text{Figure P}$$

\enddocument